\documentclass[11pt]{article}

\usepackage{graphicx,amsmath,upgreek,bm}
\usepackage{colortbl}
\usepackage[rgb, table, dvipsnames]{xcolor}
\usepackage{booktabs}
\usepackage{makecell}
\usepackage{float}
\usepackage{pifont}
\usepackage{todonotes}
\usepackage[toc,page]{appendix}
\usepackage{multirow}
\usepackage{subcaption} 	
\usepackage{amsmath,amsfonts,amsthm,bm,mathtools}
\usepackage{lipsum}
\usepackage[backend=biber,bibencoding=utf8,url=false,sorting=none,style=numeric-comp,maxbibnames = 99]{biblatex}
\usepackage{geometry}
\usepackage{algorithm}
\usepackage[noend]{algpseudocode}
\usepackage{custom_format}  %
\usepackage{arydshln}
\usepackage{quotes}
\makeatletter

\renewcommand*\env@matrix[1][*\c@MaxMatrixCols c]{%
\hskip -\arraycolsep
\let\@ifnextchar\new@ifnextchar
\array{#1}}
\makeatother
\definecolor{refcolor}{rgb}{0.65,0,0.35}
\usepackage[colorlinks=true]{hyperref}
\hypersetup{	
	linkcolor = refcolor,
	citecolor = refcolor,
	filecolor = refcolor,      
	urlcolor = blue,
	citecolor = refcolor
}

\graphicspath{ {figures/}{../figures/} }
\begin{document}
\title{\textbf{Topological Detection of Phenomenological Bifurcations with Unreliable Kernel Densities}}
\author{Sunia Tanweer \and Firas A.~Khasawneh}
\date{Dept. of Mech. Engineering, Michigan State University}

\maketitle

\abstract{
Phenomenological (P-type) bifurcations are qualitative changes in stochastic dynamical systems whereby the stationary probability density function (PDF) changes its topology. The current state of the art for detecting these bifurcations requires reliable kernel density estimates computed from an ensemble of system realizations. 
However, in several real world signals such as Big Data, only a single system realization is available---making it impossible to estimate a reliable kernel density. This study presents an approach for detecting P-type bifurcations using unreliable density estimates. The approach creates an ensemble of objects from Topological Data Analysis (TDA) called persistence diagrams from the system's sole realization and statistically analyzes the resulting set. We compare serveral methods for replicating the original persistence diagram including Gibbs point process modelling, Pairwise Interaction Point Modelling, and subsampling. We show that for the purpose of predicting a bifurcation, the simple method of subsampling exceeds the other two methods of point process modelling in performance. 
}

\section{Introduction}

Stochastic dynamics have been widely used to model systems and processes including \cite{Falk2017,Schlogl1972}, fluid dynamics \cite{VENTURI2010}, vibration control \cite{Liu2018}, energy harvesting \cite{McInnes2008, Tai2020}, genetics \cite{Lee2018}, and financial markets \cite{Chiarella2008}. Moreover, stochastic bifurcations are also seen in neurodynamical models governing our decision-making and cognitive processes~\cite{Deco2007} as well as high-dimensional neuronal models for studying conditions and onset of epilepsy~\cite{Djeundam2013, Bashkirtseva2016}. The underlying dynamics in these applications are stochastic, thus motivating qualitative and quantitative analyses of stochastic bifurcations from a variety of scientific communities \cite{Mendler2018,Arnold2013,Liu2018,Mamis2016,SchenkHoppe1996}. 

Stochastic systems can experience phenomenological bifurcation (P-type) which is characterized by a qualitative change in the topology of the joint-PDF of the response. For example, a P-bifurcation can manifest as the transition between mono-stable and bi-stable PDFs \cite{Yang2016,Yang2017,Schlogl1972,Falk2017,Arnold1996,Arnold1995}, or the emergence of stochastic limit cycles. 
These bifurcations may be induced by a change in system or excitation properties such as the damping, stiffness or intensity of noise \cite{Mamis2016, Zhao2016, Liu2020, Arnold1995,Zakharova2010,Mendler2018,Zakharova2010}. %

Since P-bifurcations are associated with qualitative changes in the PDF, their analysis has traditionally been restricted to either systems with closed form PDFs, or limited to qualitative methods that typically rely on visual inspection \cite{Kumar2022, Xu2013}. In \cite{Falk2017}, the fixed points of the PDF are computed for a one-dimensional stochastic bi-stable system. %
In \cite{Mendler2018}, the fixed points of the system PDF are used with information from limit cycles in the stochastic phase portrait to indicate ridges in the PDF. In \cite{Song2010}, a stochastic bi-stable system's bifurcations are captured in an empirical bifurcation diagram similar to those presented in \cite{Ozbudak2004}. These diagrams were used to qualitatively delineate between mono-stable and bi-stable behavior. Reference \cite{Xu2011} uses extrema of the PDF from an exact stationary solution to the stochastic Duffing-Van der Pol equation to distinguish between stochastic bi-stabililty versus monostability. In contrast \cite{Song2010,Ozbudak2004}, the analytic evaluation of \cite{Xu2011} allows the exact determination of the number of modes in the stationary distribution. %

While the methods provided in \cite{Falk2017,Mendler2018,Xu2011,Song2010,Ozbudak2004} offer some approaches for characterizing stochastic bifurcations, there are limiting factors to the current state-of-the-art. For instance, these methods only work for state spaces of at most two dimensions. If a joint-distribution (higher dimensional system) is the system's foundation, it would be much more challenging to display the shifting dynamics as is done by \cite{Song2010,Ozbudak2004}. Additionally, simply picking the number of local maxima in the PDF fails to capture richer dynamical behavior such as limit-cycles \cite{Mamis2016}. Studies which have approached bifurcations of multi-dimensional systems merely plot multiple joint distributions across a span of parameters and draw qualitative conclusions (based on visual analysis) to claim within which range of parameters a P-bifurcation occurs~\cite{Mamis2016,Djeundam2013, Liu2018}. Hence, the major hindrance in automating and standardizing these qualitative methods is their heavy reliance on visual inspection of the PDF, which makes it difficult to ascertain P-bifurcations for systems \cite{Kumar2016, Kumar2016a} with low dimensions and impossible in high dimensions. Some recent attempts at quantification have been made using changes in the sign of Shannon entropy, but the entropy-based technique is not applicable to general nonlinear dynamical systems, and bifurcation detection is particularly difficult in systems with multi-dimensional state space \cite{Venkatramani2018, Kumar2016a}. 

The authors' recent work~\cite{Tanweer2023} investigated changes in the super-level persistence diagram of the PDF or Kernel Density Estimate (KDE) of the states as the bifurcation parameter varies, and identified bifurcations through a new tool called ``homological bifurcation plot'', which utilizes the discrete difference in homology for different topological shapes to quantity a P-bifurcation. However, assumes an abundance of system realizations that enable estimating  KDE whose topology is close to the underlying PDE. However, in many real-world scenarios such as cosmological background waves or Big Data, only a single realization of the response are available. These applications are of monumental significance for progressing science, advancing healthcare systems, as well as improving national defence. For example, Big Data has the potential to revolutionize the healthcare system~\cite{Pastorino2019}, help in detecting epidemics~\cite{Aramaki2011} and improve defence and security~\cite{Xie2021}. However, the limited realizations and high dimensionality of such data lead to spurious correlations~\cite{Fan2014}, making the generation of a reliable KDE unattainable. In such cases, there is a dire need to investigate methods for detecting bifurcations by assigning statistical significance to points in the persistence diagram. %
This paper tackle this issue of tremendous significance by generating an ensemble of persistence diagrams from the one computed using the unreliable KDE, and then using statistical methods to differentiate between `signal' vs `noise'. The authors leverage prior works on modelling persistence diagrams as point processes to generate more samples~\cite{Agami2019, Adler2019, Papamarkou2022}, and those that generate more diagrams by subsampling the persistence diagram itself~\cite{Chazal2018}.  
\section{Mathematical Background}
\label{sec:maths}

This work combines multiple mathematical concepts including topological data analysis, spatial point process modelling, subsampling and outlier detection. The following subsections give a brief overview of each concept.

\subsection{Topological Data Analysis}
We start by examining one particular topological space: the cubical complex which is beneficial for discretized function (matrix) data. The subsequent subsections provide an overview of homology, superlevel-set persistent homology---directing those interested to delve into full details in works by \cite{Dey2022, Oudot2015}.

\subsubsection{Cubical Complex}
The cubical complex consists of elementary intervals: either a unit interval, denoted as $[u, u+1]$, or a degenerate interval, written as $[u,u]=[u]$ where $u$ belongs to the set of integers $\mathbb{Z}$. In an $n$-dimensional space, a cube is formed as the product $I = \prod_{i = 1}^{K}[u_i, u'_i]$, where $u'_i$ can be either $u_i$ or $u_{i+1}$. The resulting cube's dimension is determined by the count of non-degenerate intervals in this product. For instance, a $0$-cube represents a vertex, a $1$-cube an edge, a $2$-cube a square, and a $3$-cube a voxel, and so forth. Note that only coefficients in $\mathbb{Z}_2$ are used, thus disregarding the orientation of the cubes.

If one cube $\sigma$ is a subset of another $\tau$, it's stated that $\sigma$ is a face of $\tau$, denoted as $\sigma \leq \tau$. For example, both edges $\sigma_1 = [0] \times [0,1]$ and $\sigma_2 = [0,1] \times [0]$ are faces of the square $\tau = [0,1]\times [0,1]$. Consequently, a cubical complex comprises cubes where including a $d$-cube in the complex necessitates the inclusion of all its faces. 

\subsubsection{Homology}
In algebraic topology, homology serves as a mathematical tool for understanding the structure of spaces by quantifying and describing key elements like connected parts, loops, and empty regions. For our discussion, let's consider a specific fixed space denoted as $X$, typically organized as a structured pattern of simplices or cubes.

The fundamental concept behind homology revolves around examining building blocks of various dimensions—depicted as simplices or cubes—and comprehending their interconnections. These building blocks act as the fundamental constituents forming the space. The $p$th chain group, $C_p(X)$, represents combinations of $p$-dimensional building blocks, while the boundary map $\partial_p$ illustrates how these blocks link together.

Homology $H_p(X)$ is defined as the collection of cycles considering boundaries. A cycle signifies a combination of building blocks forming a closed loop or an empty space, while a boundary represents a combination that lies on the perimeter of something. Consequently, the collection of $H_p(X)$ for different values of $p$ informs us about the distinct loops, voids, or connected parts within the $X$. 

Mathematically, this concept is expressed as $H_p(X) = \Ker(\partial_p) / \im(\partial_{p+1})$, where $\Ker(\partial_p)$ denotes the set of cycles (combinations forming closed structures), and $\im(\partial_{p+1})$ represents the set of boundaries (combinations situated on the edge). Thus, the quotient, $H_p(X)$, captures the unique and non-duplicative aspects of the space at the $p$-dimensional level. Informally, the homology of a space provides a vector describing topological characteristics inherent to the space with each homology dimension $p$ providing some aspect of the shape of the underlying space such as connected components ($p=0$), loops ($p=1$), voids ($p=2$), or higher-dimensional counterparts for $p>2$. 

\subsubsection{Superlevel Persistent Homology}
Persistent homology is a modern variant of homology that is a widely employed technique in Topological Data Analysis (TDA), captures details about the structure of a parameterized space by tracing the changes in its homology as the parameter varies \cite{Dey2022}.

Consider a topological space equipped with a function $f:X \to \R$. Given data comprised of a series of real values $a_1\leq a_2 \leq \cdots \leq a_n$, we can define two filtrations: the \textit{sub-level-set} and the \textit{super-level-set}. In the former, $X_a = f^{-1}(-\infty,a]$ forms a sequence $X_{a_1} \subseteq X_{a_2} \subseteq \cdots \subseteq X_{a_n}$. Conversely, the latter creates $X^a = f^{-1}[a,\infty)$, forming $X^{a_n} \subseteq X^{a_{n-1}} \subseteq \cdots \subseteq X^{a_1}$.

A fundamental theorem in algebraic topology asserts the functoriality of homology. In the context of inclusions $X \subset Y$, it implies an induced map (a linear map on vector spaces) on relevant homology groups $H_p(X) \to H_p(Y)$. Leveraging this, we establish sub-level-set and super-level-set homology filtrations, progressing from $H_p(X_{a_1})$ to $H_p(X_{a_n})$ and $H_p(X^{a_n})$ to $H_p(X^{a_1})$, respectively. 

\begin{figure}[!htbp]
    \centering
    \includegraphics[height = 1.5in]{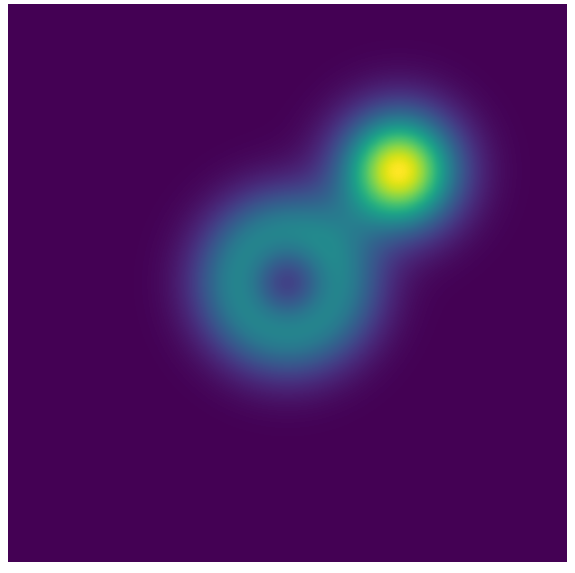} \includegraphics[height = 1.5in]{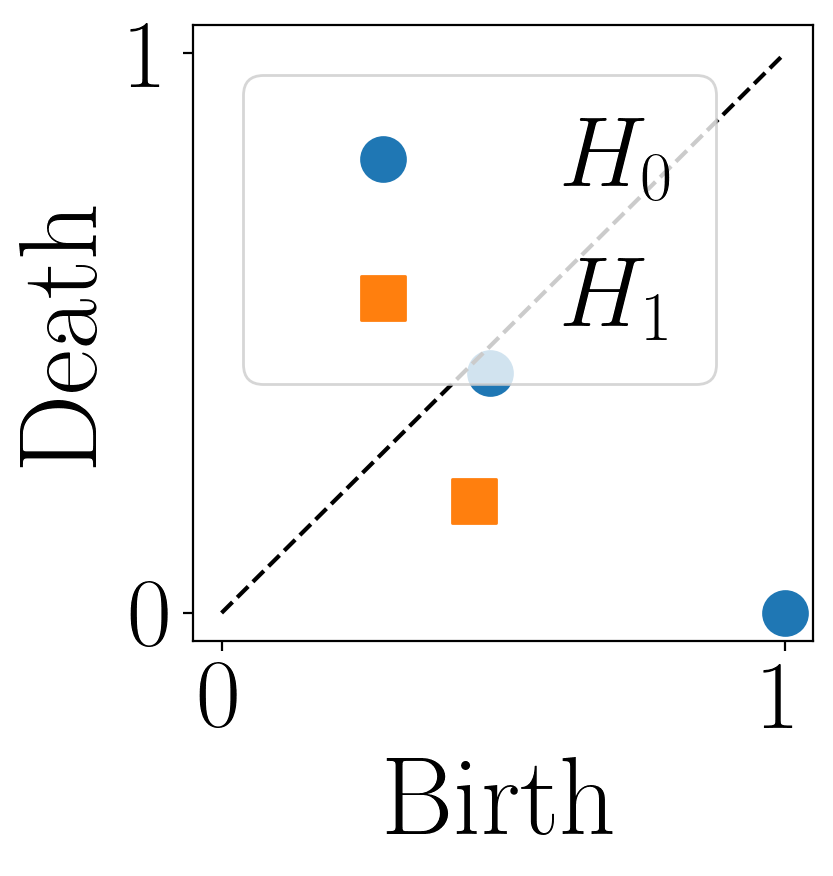}\\
    \includegraphics[height = 0.75in]{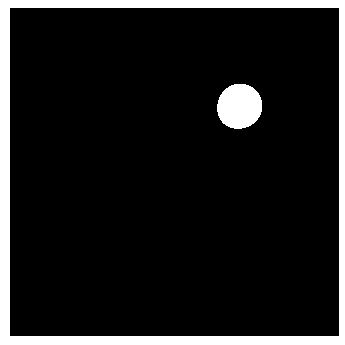}
    \includegraphics[height = 0.75in]{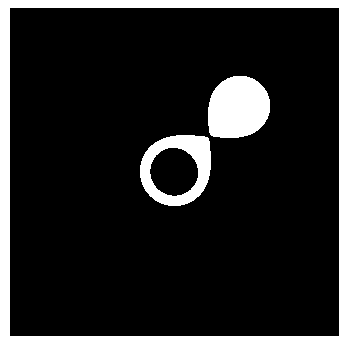}
    \includegraphics[height = 0.75in]{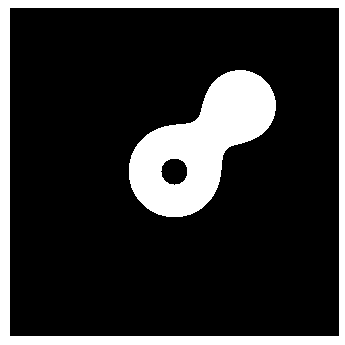}
    \includegraphics[height = 0.75in]{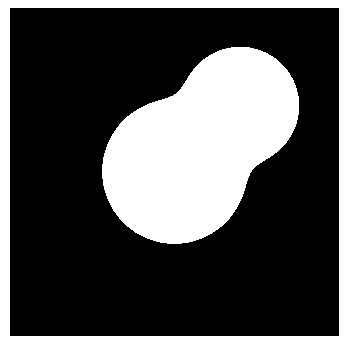}
    \caption{Superlevel cubical persistence of unit-normalized image data shown. Figures correspond to (a) $K^{0.8}$, (b) $K^{0.5}$, (c) $K^{0.4}$ and (d) $K^{0.2}$.}
    \label{fig:CubicalPersistence}
\end{figure}

Persistent homology further stipulates that any filtration of this form can be uniquely decomposed into pairs $(b,d)$. Here, $b=a_i$ denotes the emergence of a new homology class in $H_p(X_{a_i})$ or $H_p(X^{a_i})$, distinct from the image of $H_p(X_{a_{i-1}})$ or $H_p(X^{a_{i+1}})$. Conversely, $d=a_j$ signifies the merging of a class with the image of an older class entering $H_p(X_{a_j})$ or $H_p(X^{a_j})$. 

This data is typically visualized in a persistence diagram, depicting pairs as points in a plane. Notably, in super-level-set persistence, $b>d$, where $b$ indicates the parameter's value when a structure first appears (termed "birth"), and $d$ denotes its disappearance ("death"). Consequently, super-level-set points reside below the diagonal line $\Delta = {(x,x) \mid x \in \R}$.

In this paper, we utilize particularly one persistence: the super-level-set Persistence of Cubical Complexes. For this, consider a grayscale image or an $m \times n$ matrix $M$, defining its domain as the cubical set $K=\mathcal{K}([0,m]\times[0,n])$. Here, the image is represented as a function $f:K \to \mathbb{R}$ on this cubical set, where for a given cube $s_{i,j}=[i,i+1]\times[j,j+1]$ in $K$, $f(s_{i,j})$ corresponds to $M_{i,j}$, the matrix entry value. The super-level-set of $f$, denoted as $K^a = f^{-1}[a, \infty)$, also forms a cubical complex. Analyzing the super-level-set persistence of this function allows us to observe its evolution as the parameter (here, $a$) changes, yielding different cubical complexes and their associated topological features. Fig.~\ref{fig:CubicalPersistence} shows an example of an image data with its evolving homology. At $a \sim 0.8$, we see that one $H_0$ component exists. At $a \sim 0.5$, we see that another $H_0$ component takes birth along with an $H_1$ component. At $a \sim 0.4$, we notice the two $H_0$ components merge into one leading to the death of one of them. Finally, at $a \sim 0.2$, the loop in the image fills up causing death of the $H_1$ component. 

\subsubsection{Detection of Bifurcation using Superlevel Persistence}

\begin{figure}[!htbp]
    \centering
    \includegraphics[width = 0.9\linewidth]{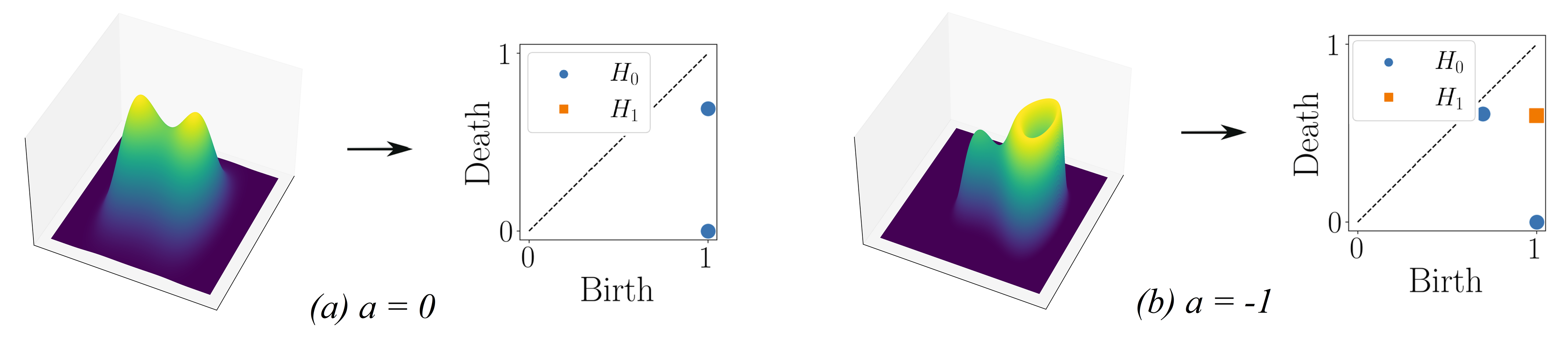}
    \caption{Two PDFs with their cubical persistence diagrams. (a) A bistable PDF has two $H_0$ components in the persistence diagram and no $H_1$ component while (b) a monostable PDF with limit cycle oscillations will have two $H_0$ components for each peak and one $H_1$ component for the limit cycle.}
    \label{fig:CubPerBif}
\end{figure}

In~\cite{Tanweer2023} the authors showed the correspondence between the number of peaks in a reliable PDF and the number of $H_0$ components, and the number of limit cycle oscillations and the number of $H_1$ components. That correspondence allowed contrasting different persistence diagrams associated to topologically different reliable PDFs to detect P-bifurcations, see Fig.~\ref{fig:CubPerBif} for an example. 

\subsection{Spatial Point Process Modelling}
A spatial point pattern $\chi^{(0)} \in \mathbb{R}^d$ can be modelled as one realization of a random process called a spatial point process. One of the primary goals of this model is to generate more realizations ($\chi^{(1)}, \chi^{(2)}, \cdots, \chi^{(n}$) of point patterns. The modeling of the underlying random process depends on several factors, including the inherent characteristics of the data and the nature of the phenomena being studied. As such, many models have been proposed to this day like the Poisson Point Process, the Binomial Point Process or the Gibbs Point Process, more details can be found in~\cite{Baddeley2012, Baddeley2005}. This paper adopts modelling persistence diagrams as point processes to generate more realizations. We use two point processes that we describe below: Gibbs (GPP) and Pairwise Interaction Point Process (PIPP). 

\subsubsection{Gibbs Point Process (GPP)}
\begin{figure}[!htbp]
    \centering
    \includegraphics[width = 0.95\linewidth]{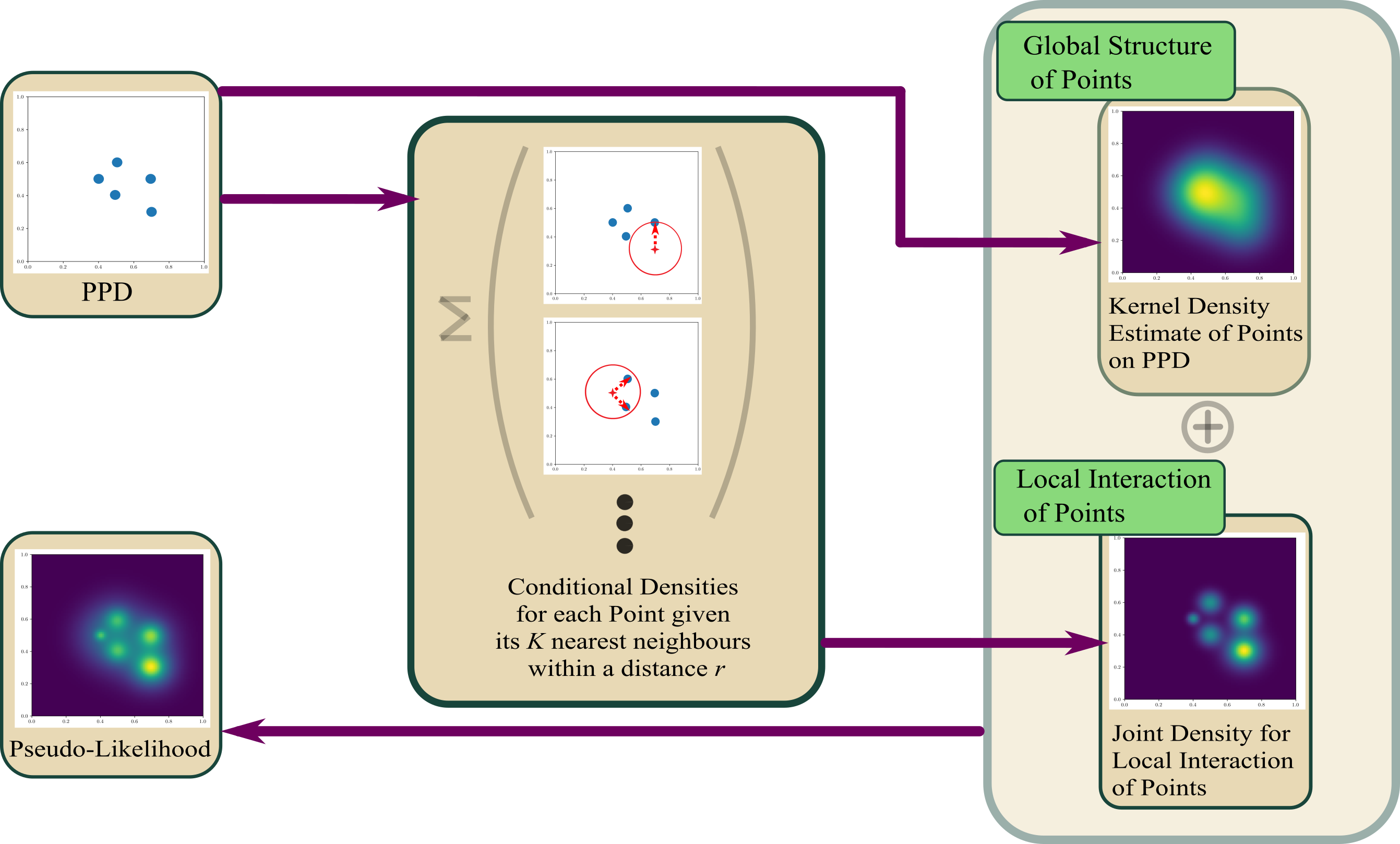}
    \caption{Gibbs point process modelling. Given a PPD, the local interaction between its points is calculated as a conditional density for each point given its neighbors and the global interaction is calculated as a kernel density. The two interactions combined lead to a pseudo-likelihood estimate for the PPD.}
    \label{fig:gibbs_methods}
\end{figure}

The Gibbs process for modelling persistence diagrams was first introduced in~\cite{Adler2017} and later refined in~\cite{Adler2019}. The modelling method is summarized in Fig.~\ref{fig:gibbs_methods} where given a projected persistence diagram (PPD) with birth-lifetime coordinates a Gibbs model computes a local and a global generalization for it. The local terms capture the interaction of each point with its neighbours within a neighborhood while the global term is essentially a density function estimated on top of the PPD. Both these terms are used to estimate a pseudo-likelihood which can then be utilized for proposing new points and determining their acceptance probability. %
Once the pseudo-likelihood is computed,~\cite{Adler2019} uses a Metropolis-Hastings Markov Chain Monte Carlo (MH-MCMC) approach for sampling new PPDs. This only leads to relocation of randomly chosen points from the PPD based on the acceptance probability without increasing or decreasing the points in the PPD.

\subsubsection{Pairwise Interaction Point Process (PIPP)}
\begin{figure}[!htbp]
    \centering
    \includegraphics[width = 0.99\linewidth]{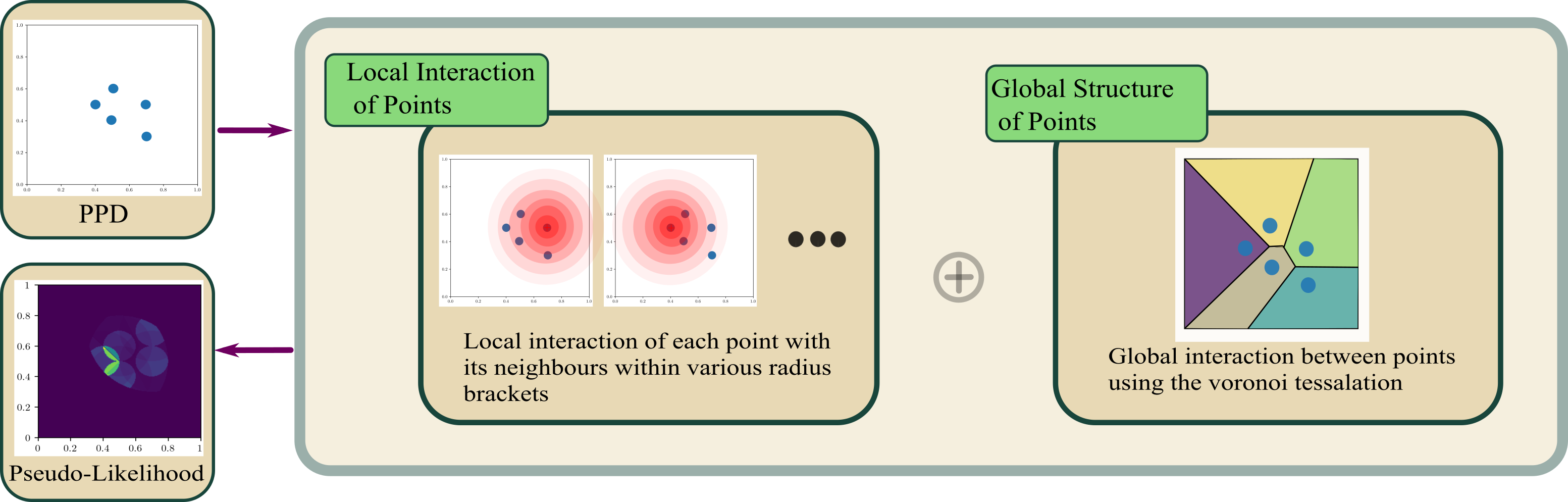}
    \caption{PIPP modelling for persistence diagrams. Given a PPD, the local interaction between points is estimated using its neighbours within different radii brackets while the global interaction is estimated using the areas of voronoi tessalations on the PPD. Combining these interactions gives a pseudo-likelihood estimate for the PPD.}
    \label{fig:pipp_method}
\end{figure}

The work in~\cite{Papamarkou2022} builds upon~\cite{Adler2019} and uses Pairwise Interaction Point Process (PIPP) to decompose a given PPD into its local and global interactions. Figure~\ref{fig:pipp_method} showcases a summary of PIPP modelling where the main difference between the models in \cite{Papamarkou2022} and \cite{Adler2019} lies in how the local and global terms are obtained. For the global part, PIPP uses the areas from a voronoi tessalation implementation on the given PPD, whereas the local part is modelled using a piecewise-constant pairwise interaction function which places multiple circles of increasing radii on a specific point in the PPD and compute its interaction with other points based on that radial distance. Then the global and local terms are merged into a pseudo-likelihood similar to the Gibbs case. 
Once the pseudo-likelihood is computed,~\cite{Papamarkou2022} uses a Reversible Jump Markov Chain Monte Carlo (RJ-MCMC) approach for sampling new PPDs. This not only relocates randomly chosen points from the PPD based on the acceptance probability but may also add or subtract points from the original PPD in subsequent replicates.

\subsection{Subsampling the Persistence Diagram}
\begin{figure}[!htbp]
    \centering
    \includegraphics[width = 0.85\linewidth]{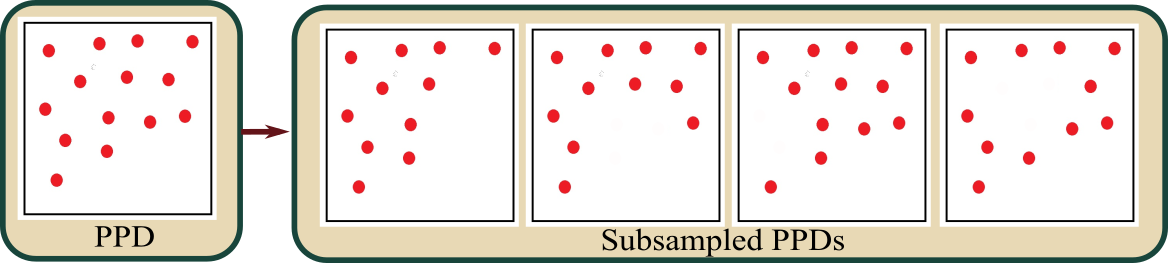} 
    \caption{Given a PPD, random subsampling of points with replacement can be used to generate further realizations of the PPD. All of these new samples can be used to generate statistically meaningful results.}
    \label{fig:bootstrap_method}
\end{figure}

Figure~\ref{fig:bootstrap_method} depicts subsampling of PPDs, which was introduced by Chazal et al.~\cite{Chazal2018}. Subsampling implements a bootstrapping methodology to supplement the single PPD realization for generating more realizations. Using random subsampling with replacement many new PPDs can be generated from one PPD and used to obain statistically meaningful results for the underlying unreliable KDE. 

\subsection{Noise Separation by Outlier Detection}
Informative points in a PPD typically lie away from the birth-axis and are small in number while the unreliable points associated with noise often tend to be close to the birth-axis and are numerous. Consequently, outlier detection methods can be used for defining a threshold that separates outliers which can be considered or important from the noise below the threshold. We use two methods of outlier detection with both GPP and PIPP models: bootstrapping of lifetimes in the persistence diagrams, and outlier detection using Mahalanobis distance (MD). The former is a widely used approach, so we only focus on describing the latter. 

MD is an effective multivariate distance metric that measures the distance between a point $x$ and a distribution according to
\[MD = \sqrt{(x-m)^T C^{-1} (x - m)}\]
where $m$ is the mean value of the independent variable, and $C$ is the covariance matrix. Once the distance is computed for each point in a certain PPD, the outliers can be found as points with a high MD (typically $3$ standard deviations away from the average MD of the dataset).

Since the points in the PPD associated to noise tend to be many and concentrated near the birth-axis, they are expected to have a high density which overshadows that of informative points in the PPD. Therefore, we use MD to separate the two types of points in the PPD.

\section{Methods}

\begin{figure}[!htbp]
    \centering
    \includegraphics[width=0.7\linewidth]{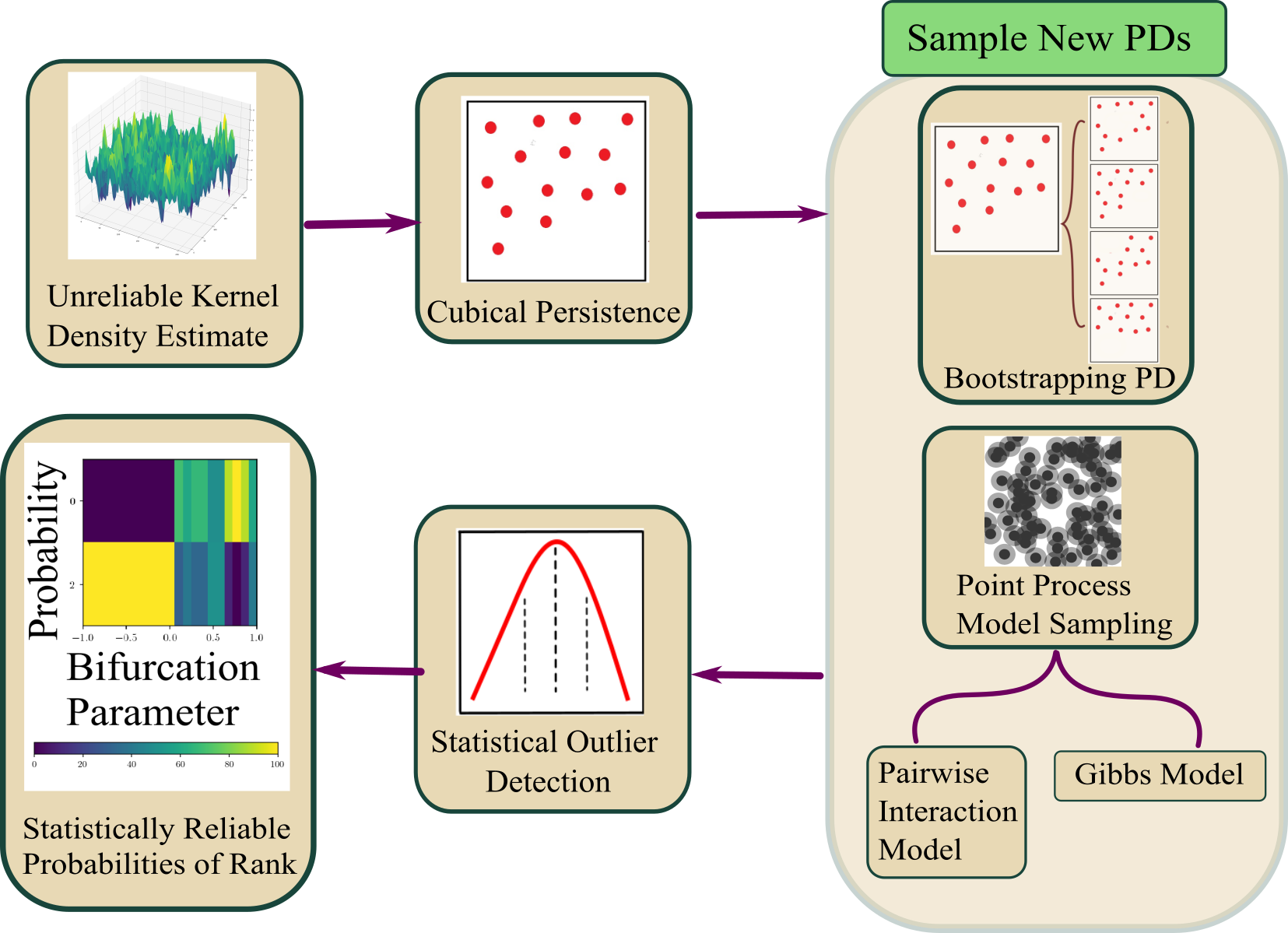}
    \caption{Overview of the method for bifurcation detection given unreliable kernel density estimates.}
    \label{fig:method}
\end{figure}

Figure~\ref{fig:method} gives an overview of our method for detecting bifurcations given an unreliable density estimate given a single realization $X$. Once the KDE has been estimated it is unit normalized, and the corresponding superlevel persistence diagram $D \in [0, 1]^2$ is computed. 

Next we generate more persistence diagrams given the diagram $D$. This can be achieved using any of the three methods discussed in Section~\ref{sec:maths}, such as Gibbs or PIPP modelling, or subsampling. Once an ensemble of persistence diagrams is obtained, any of the outlier detection methods that we described in Section~\ref{sec:maths} can be used to determine the significant points in each PPD. %
This method is repeated for each system realization over the range of bifurcation parameters to build a probability versus homology rank plot which can inform us of the occurrence of a bifurcation based on significant changes in the probability of a certain homology rank. 

\section{Results and Discussion}
\label{sec:results}

We use all three methods of persistence diagram replication given unreliable KDEs from single system realizations to study P-bifurcations in stochastic Duffing, Rayleigh Vander-Pol and Quintic oscillators. From each original persistence diagram, $500$ further realizations were generated. This section compares the results from all these methods (3 for replication, 2 for outlier detection) for each of these. Since our inputs are unreliable the it is impossible to report the exact parameter value where the bifurcation occurred. Instead, we compare the three methods to see which one gives a clearer indication of a change in the system behaviour close to our bifurcation point. 

\subsection{Stochastic Duffing Oscillator}

The stochastic Duffing oscillator is represented by 
\[\ddot{X} + \dot{X} + hX + X^{3} = q_{1}dW_{1}.\]

where $q_1$ is the amplitude of the added Gaussian white noise and $h$ is the bifurcation parameter such that for $h < 0$ the system has a monostable PDF, and for $h \geq 0$ the system has a bistable PDF. In~\cite{Tanweer2023}, it was shown that for $h < 0$, the persistence diagram had $1$ $H_0$ component and for $h \geq 0$, the persistence diagram had $2$ $H_0$ components. 

\begin{figure}[!htbp]
\centering
\includegraphics[width=0.3\linewidth]{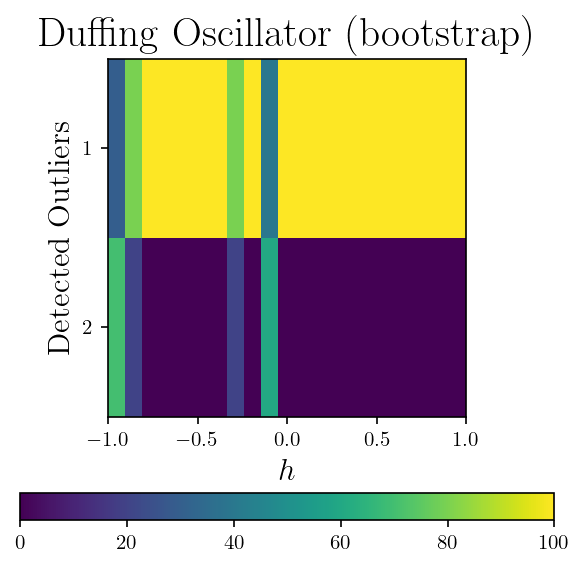}
\includegraphics[width=0.3\linewidth]{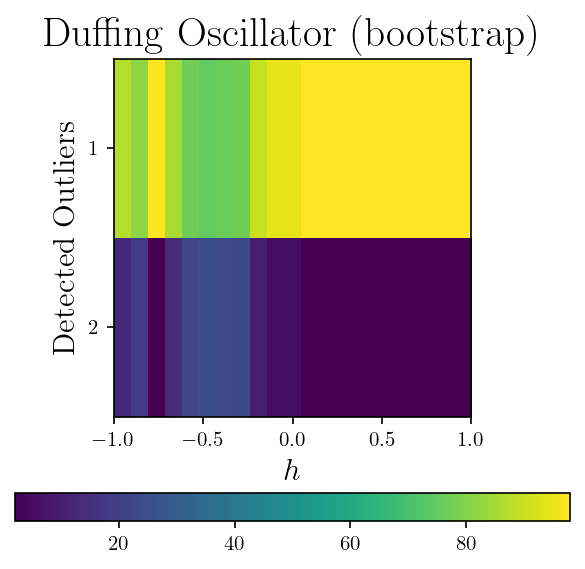}
\includegraphics[width=0.3\linewidth]{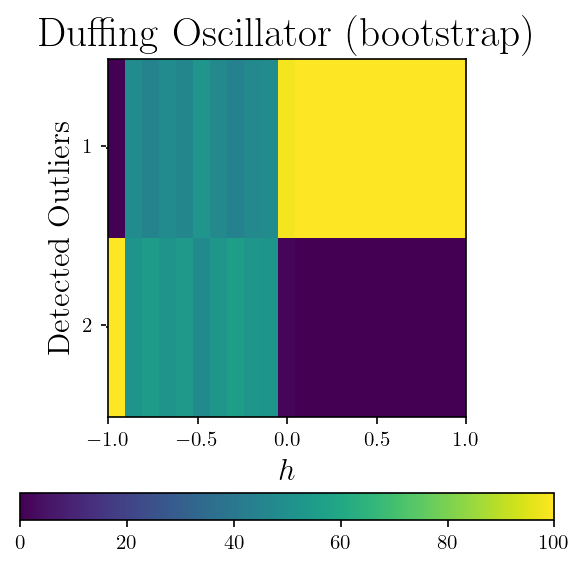}
\\
\includegraphics[width=0.3\linewidth]{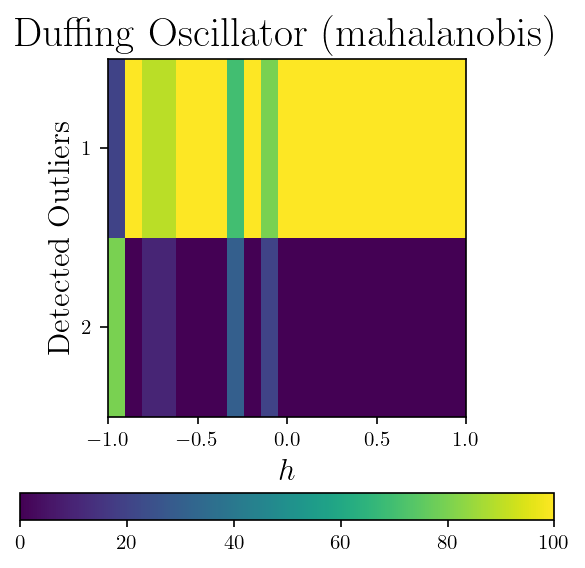}
\includegraphics[width=0.3\linewidth]{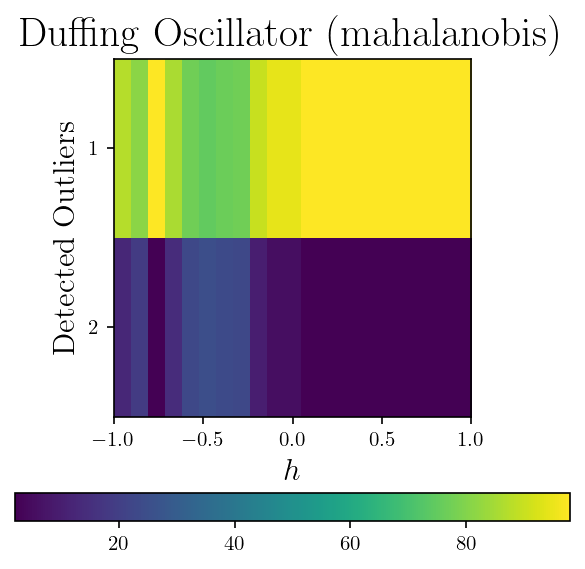}
\includegraphics[width=0.3\linewidth]{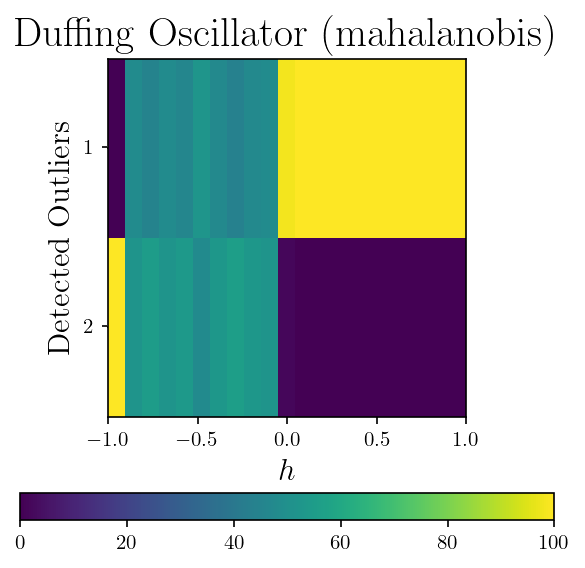}
\caption{Comparison of $H_0$ plots for Duffing Oscillator for outlier detection. Left to Right: Gibbs, PIPP and Subsampling.}
\label{fig:duffing}
\end{figure}

Figure~\ref{fig:duffing} shows the probabilistic plots for $H_0$ classes from all $500$ replicates generated by the three methods (GPP, PIPP, subsampling). Regardless of the replication method, all instances detect a monostability for $h \geq 0$ as evident by the $100\%$ probability for $1$ outlier and $0\%$ probability for $2$ outliers. For $h < 0$, we expect a bistable system and wish to see an increase in the probability for $2$ outliers. Figure~\ref{fig:duffing} shows that for $h < 0$, the probability for homology of $2$ slowly increases for the subsampling results, followed by PIPP replication, indicating an increased probability of a bifurcation. In contrast, the GPP modelling results are unable able to generate replicates with $2$ outliers for the majority of $h<0$ regime regardless of the outlier detection method. 

\subsection{Stochastic Rayleigh Vander-Pol Oscillator}

The stochastic Rayleigh Vander-Pol oscillator is given by 
\begin{equation*}
  \ddot{X} + (h + X^{2} + \dot{X}^{2})\dot{X} + X = q_{1}dW_{1},
\end{equation*}
where $q_1$ and $h$ are the same as those defined for Duffing earlier. In this case, the system is monostable for $h \geq 0$ and exhibits a stochastic limit cycle for $h < 0$. In~\cite{Tanweer2023}, it was shown that for $h < 0$, the persistence diagram had $1$ $H_0$ and $1$ $H_1$ component (corresponding to the limit cycle) and for $h \geq 0$, the persistence diagram had $1$ $H_0$ component with no $H_1$ components. 

\begin{figure}[!htbp]
\centering
\includegraphics[width=0.3\linewidth]{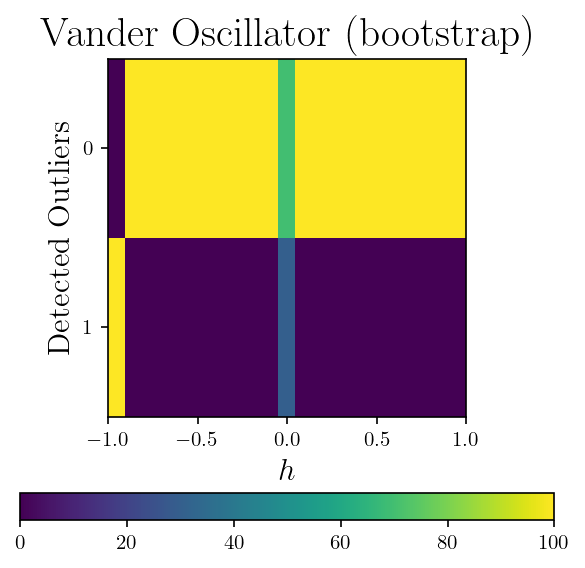}
\includegraphics[width=0.3\linewidth]{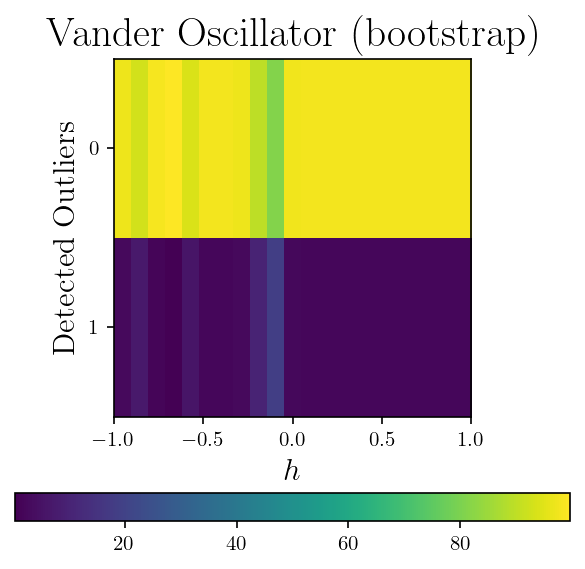}
\includegraphics[width=0.3\linewidth]{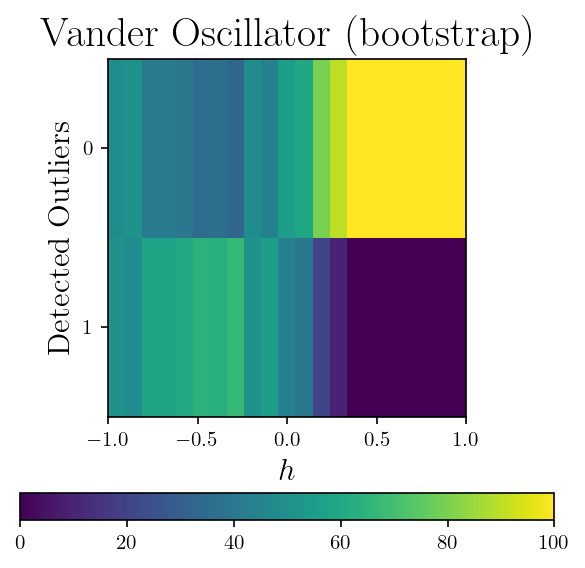}
\\
\includegraphics[width=0.3\linewidth]{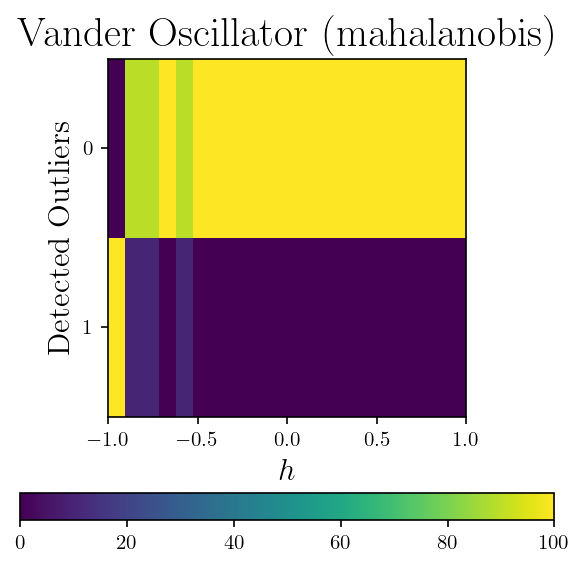}
\includegraphics[width=0.3\linewidth]{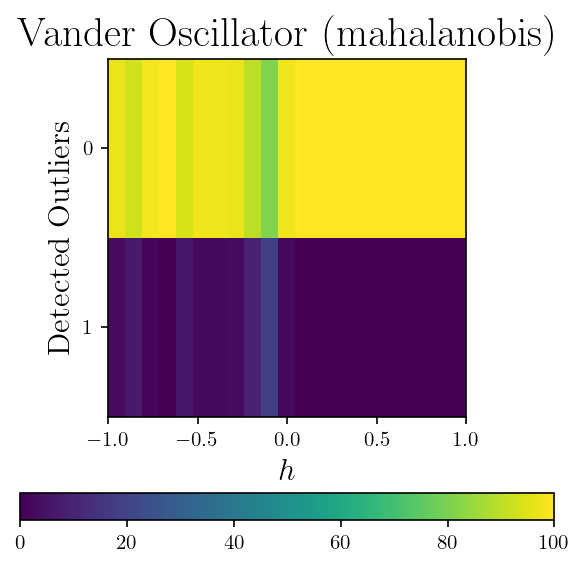}
\includegraphics[width=0.3\linewidth]{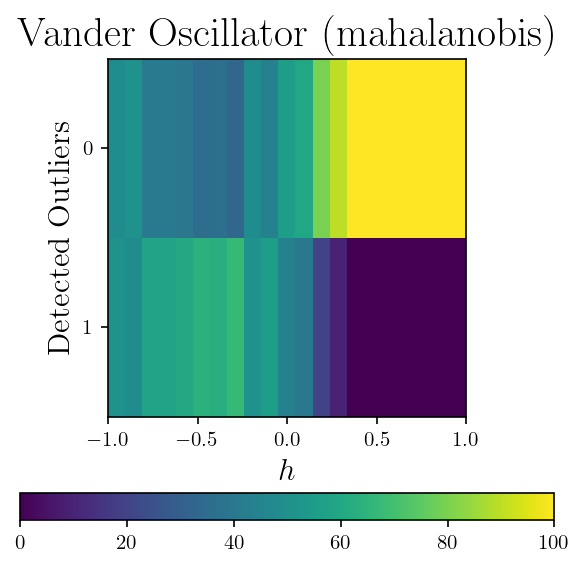}
\caption{Comparison of $H_1$ plots for Vander Oscillator. Left to Right: Gibbs, PIPP and Subsampling.}
\label{fig:vander}
\end{figure}

Fig.~\ref{fig:vander} shows the probabilistic plots for $H_0$ classes from all $500$ replicates generated by the three methods. Regardless of the replication method, all instances detect no limit cycle for $h \geq 0$ as evident by the $100\%$ probability for $1$ outlier and $0\%$ probability for $2$ outliers.
Again, the results by GPP are not capable of detecting the bifurcation accurately. The results from PIPP are somewhat better and depict slight changes in the probability for $0$ and $1$ rank of homology as the $h<0$ threshold is crossed but the probability changes are insignificant. 
In contrast to GPP and PIPP, the method of Subsampling shows clear and significant changes in the probability for $h<0$, depicting a bifurcation---even though this change starts manifesting at $h \sim 0.2$ instead of $h = 0$ (which is the true bifurcation point). 

\subsection{Stochastic Quintic Oscillator}

The stochastic Quintic oscillator is forced by both additive and multiplicative Gaussian noises and is given by
\begin{equation*}
\begin{split}
  \ddot{X} + 2[D_{11}(2U(X) + h) - \frac{1}{2}D_{22}]\dot{X} + \\
  2[D_{22}(2U(X) + h) + D_{11}]\dot{X}^{3} + 2D_{22}\dot{X}^{3} + h_{0}(X) \\
  = dW_{1} + \dot{X}dW_{2},
\end{split}    
\end{equation*}
where $h_{0}(X) = x^{3}_{1} + ax^{2}_{1} - x_{1}$. 

For this oscillator, both $h$ and $a$ can be treated as bifurcation parameters keeping the other constant. Two cases were considered in~\cite{Tanweer2023}:
\begin{enumerate}
    \item Varying $h \in [-1, 1]$ with $a=0$,
    \item Varying $a \in [-1, 0]$ with $h=1$.
\end{enumerate}

\subsubsection{Varying \texorpdfstring{$h \in [-1, 1]$ with $a=0$}{ h in [-1,1] with a=0}}
For bifurcation parameter $h$, the system PDF evolves from a bistability at $h=1$ to an inverted bistability at $h=-1$. In~\cite{Tanweer2023}, we see the persistence diagram evolving from $2$ $H_0$ and no $H_1$ points to $1$ $H_0$ and $2$ $H_1$ points respectively. 

\begin{figure}[!htbp]
\centering
\includegraphics[width=0.3\linewidth]{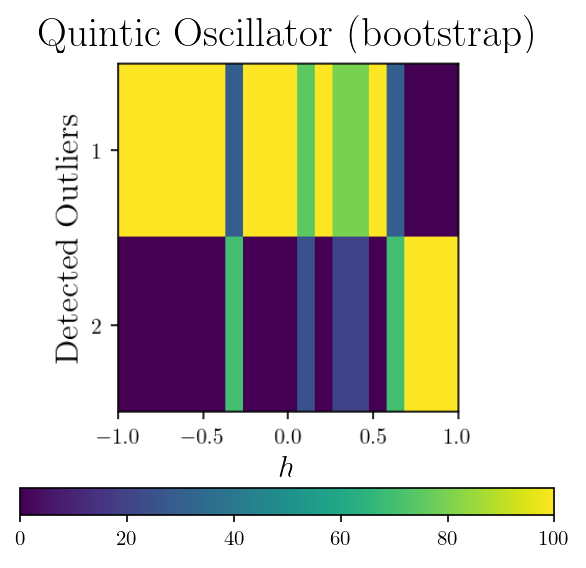}
\includegraphics[width=0.3\linewidth]{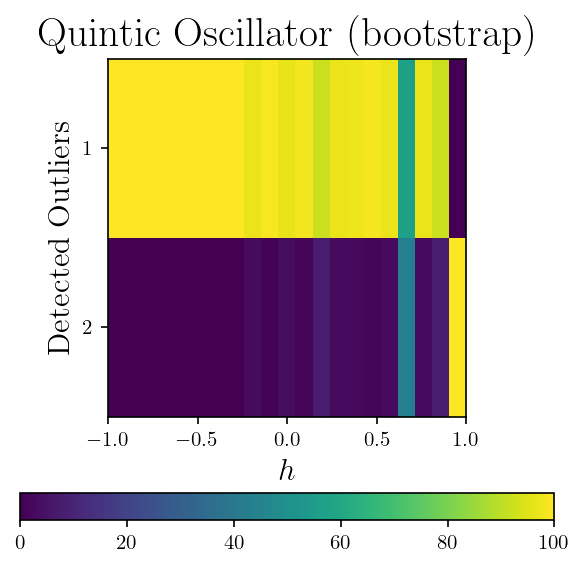}
\includegraphics[width=0.3\linewidth]{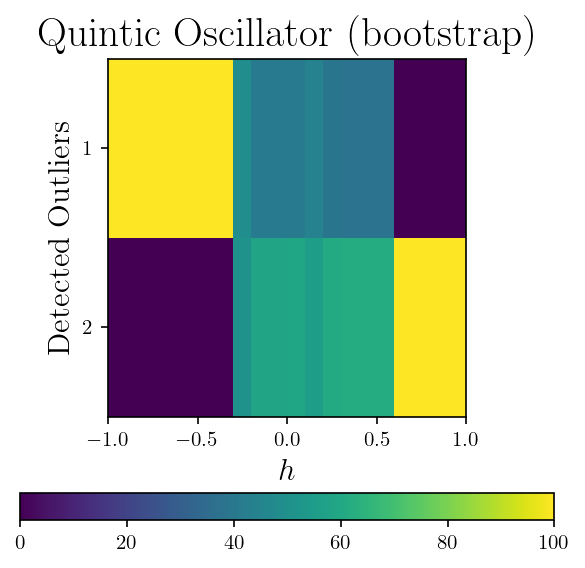}
\\
\includegraphics[width=0.3\linewidth]{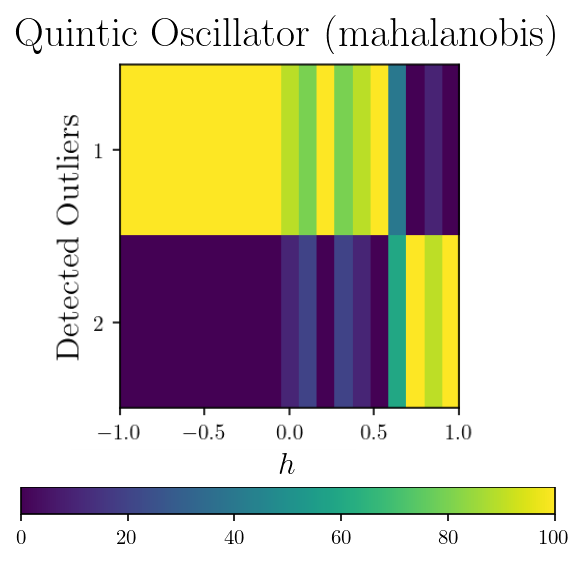}
\includegraphics[width=0.3\linewidth]{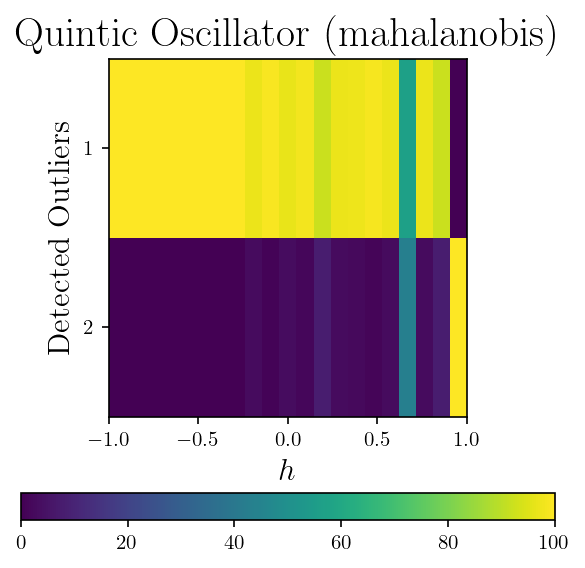}
\includegraphics[width=0.3\linewidth]{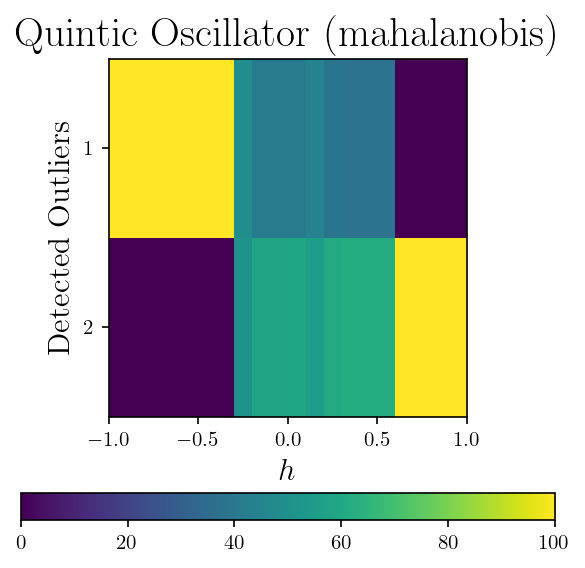}
\caption{Comparison of $H_0$ plots for Quintic Oscillator for $a = 0$ and $h \in [-1, 1]$. Left to Right: Gibbs, PIPP and Subsampling.}
\label{fig:quintic_a_H0}
\end{figure}

Fig.~\ref{fig:quintic_a_H0} shows the probabilistic plots for $H_0$ classes from all $500$ replicates generated by the three methods (GPP, PIPP, subsampling). Contrary to the previous cases, both GPP and PIPP qualitatively detect a change in the system's PDF close to $h=0.5$ (earlier than the actual bifurcation point of $h=0$), whereas the subsampling method gives a band of $h \in [-0.3, 0.5]$ with probabilities other than $0\%$ or $100\%$, indicating the occurrence of a bifurcation somewhere inside this band. 

\begin{figure}[!htbp]
\centering
\includegraphics[width=0.3\linewidth]{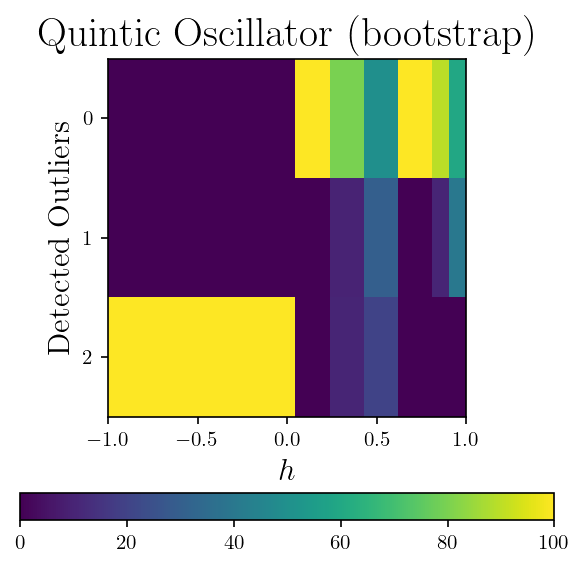}
\includegraphics[width=0.3\linewidth]{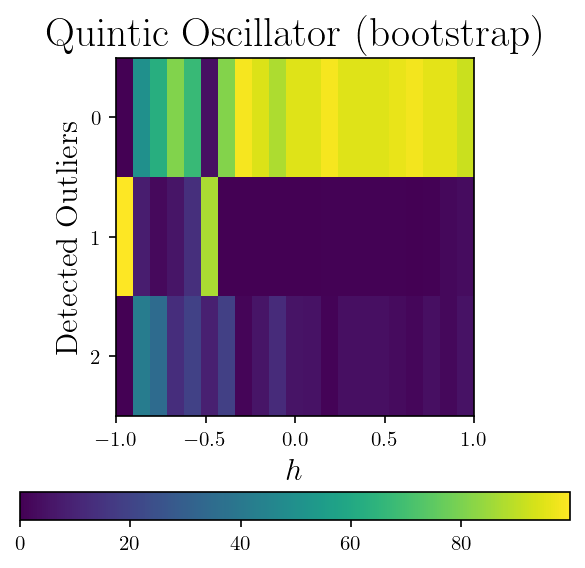}
\includegraphics[width=0.3\linewidth]{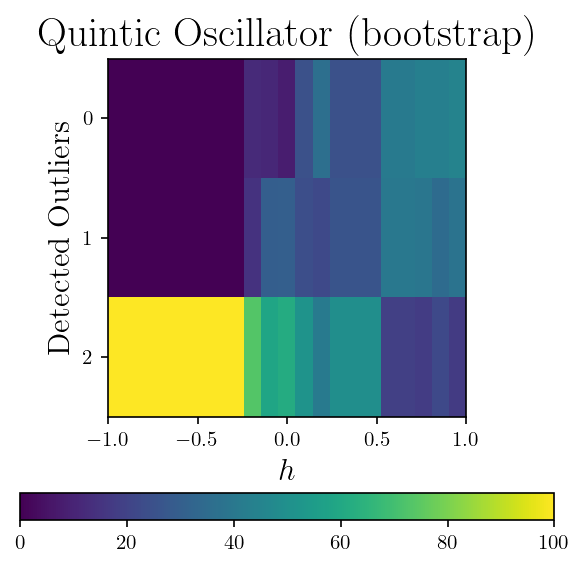}
\\
\includegraphics[width=0.3\linewidth]{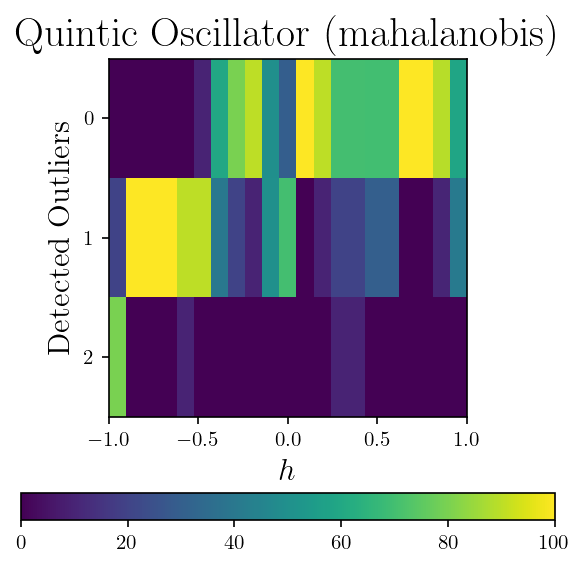}
\includegraphics[width=0.3\linewidth]{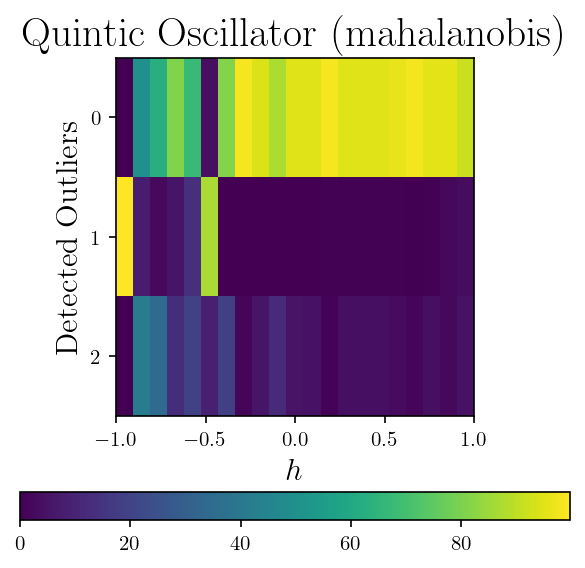}
\includegraphics[width=0.3\linewidth]{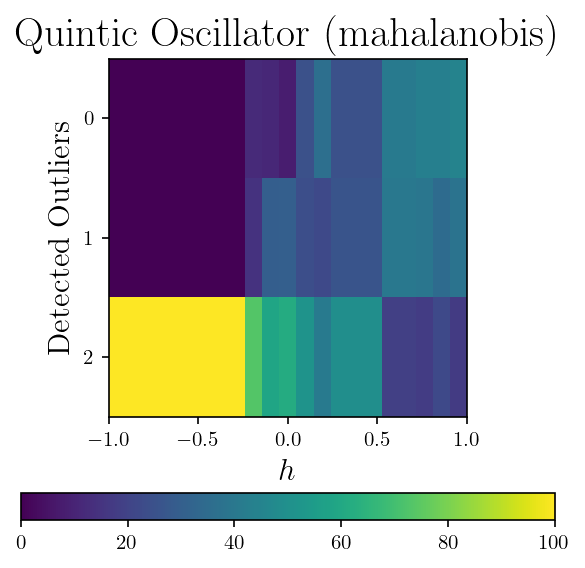}
\caption{Comparison of $H_1$ plots for Quintic Oscillator for $a = 0$ and $h \in [-1, 1]$. Left to Right: Gibbs, PIPP and Subsampling.}
\label{fig:quintic_a_H1}
\end{figure}

Figure~\ref{fig:quintic_a_H1} shows the probabilistic plots for $H_1$ classes from all $500$ replicates generated by the three methods (GPP, PIPP, subsampling). For this original PD of $H_1$ class, GPP with bootstrapping outlier detection performs much better than PIPP by detecting a clear change at $h=0$. Subsampling, like before, prematurely starts detecting the bifurcation slightly around $h \sim -0.2$ rather than at the actual bifurcation point $h=0$. 

\subsubsection{Varying \texorpdfstring{$a \in [-1, 0]$ with $h=1$}{a in [-1, 0] with h=1}}
For the bifurcation parameter $a$, the system's PDF transforms from a bistability at $a=0$ to a monostability with limit cycle at $a=-1$, while the persistence diagram evolves from having $2$ $H_0$ points to having $2$ $H_0$ and $1$ $H_1$ point respectively. 

\begin{figure}[!htbp]
\centering
\includegraphics[width=0.3\linewidth]{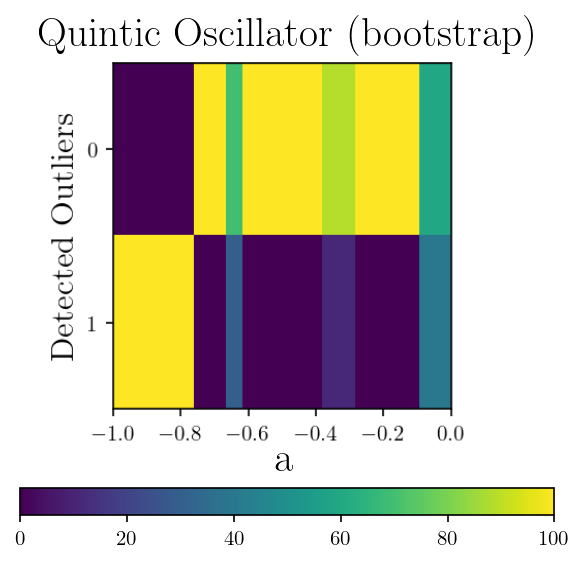}
\includegraphics[width=0.3\linewidth]{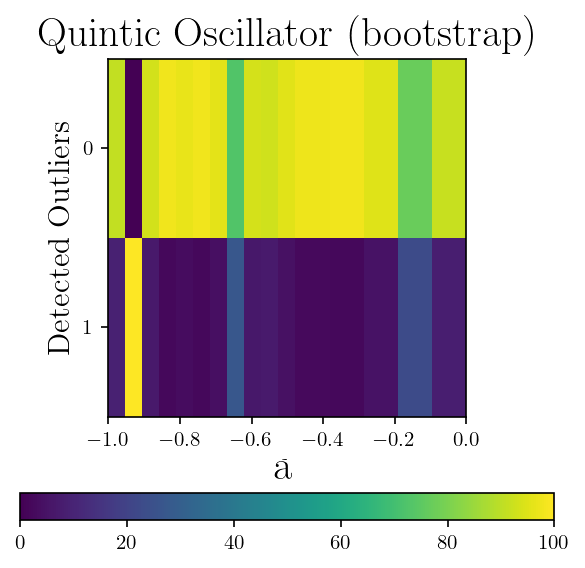}
\includegraphics[width=0.3\linewidth]{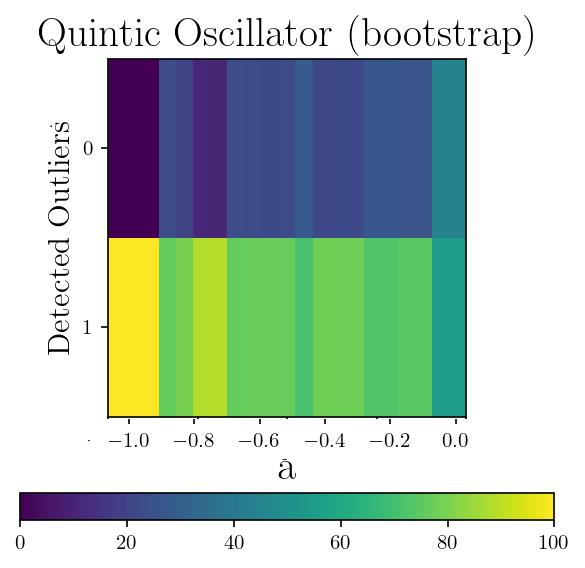}
\\
\includegraphics[width=0.3\linewidth]{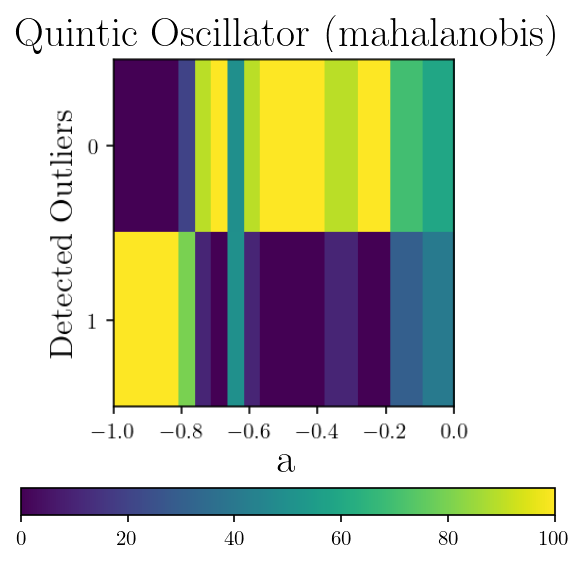}
\includegraphics[width=0.3\linewidth]{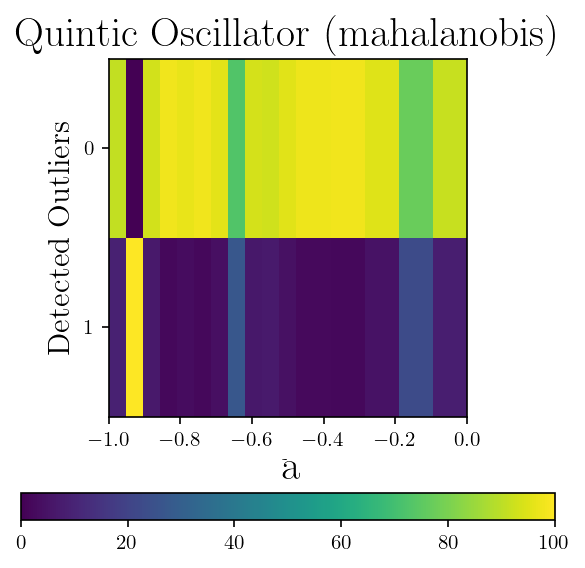}
\includegraphics[width=0.3\linewidth]{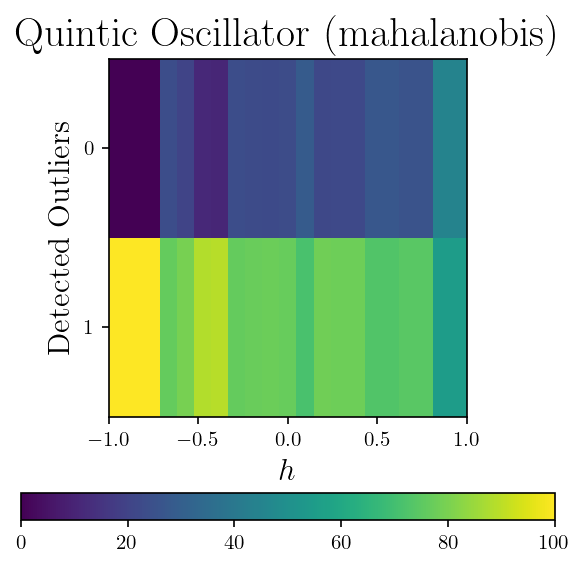}
\caption{Comparison of $H_1$ plots for Quintic Oscillator for $h = 1$ and $a \in [-1, 0]$. Left to Right: Gibbs, PIPP and Subsampling.}
\label{fig:quintic_h_H1}
\end{figure}

Figure~\ref{fig:quintic_h_H1} shows the probabilistic plots for $H_1$ class from all $500$ replicates generated by the three methods (GPP, PIPP, subsampling). For this original PD of $H_1$ class, both GPP and PIPP do not detect the bifurcation close to $a \sim -0.6$ very well. Subsampling, like earlier cases, begins showing a change in probabilities close to the bifurcation point, albeit a little off from the actual point.

In each of these scenarios of specific to foretelling a bifurcation, Subsampling has consistently demonstrated more reliable performance compared to GPP and PIPP. Notably, it provides a smoother variation in probabilities, contrasting with the occasionally unpredictable results of GPP and PIPP. 
One possible explanation is that the PDs from our oscillators have very few true points (usually $0$, $1$ or $2$) compared to the many noise points lying in the same or nearby region on the PD space. Since both GPP and PIPP are based on generating a likelihood for the original PD, the huge difference in the densities of the true and noise points---in many cases---causes the true points to be ignored in the replicates generated due to their low likelihood of occurrence. However, in the case of Subsampling, a true point is as likely to be removed in any new iteration as any other noise point, making the subsampling approach much better for our oscillators where the true-to-noise ratio is much low.

\section{Conclusion}
\label{sec:conc}

We investigated the challenging task of P-bifurcation detection in scenarios where reliable kernel density estimates (KDE) are infeasible due to the scarcity of realizations. The comparison of three distinct methods---Gibbs point process modeling (GPP), pairwise interaction point process modeling (PIPP) and subsampling---demonstrated their efficacy in replicating the persistence diagram (PD) from the unreliable KDEs. Notably, subsampling consistently yielded the closest indication of a P-bifurcation among these methodologies, showcasing its robustness in detecting a bifurcation albeit a little earlier than the exact point of bifurcation. However, under the heavy constraint of unreliable KDE from single system realization, probabilistic bands such as the ones seen in the subsampling case can provide actionable information for decision making to counter an imminent bifurcation. Moreover, the study incorporated two outlier detection techniques, specifically bootstrapping and Mahalanobis distance. Their comparable performances in most cases underline their utility in identifying outliers within this context, irrespective of the chosen method for PD replication. 

While the outcomes showcased a consistent advantage for subsampling, the variable performance of GPP and PIPP warrants further investigation. The intermittently better results of GPP or PIPP suggest potential nuances in their applicability based on specific data characteristics or system behaviors. Overall, this research not only sheds light on the intricate landscape of bifurcation detection in the presence of unreliable KDEs but also emphasizes the robustness well-established methodologies such as subsampling. 
Future endeavors could delve deeper into understanding the idiosyncrasies influencing the performance variations of GPP and PIPP, paving the way for refined methodologies tailored to diverse stochastic systems in fields such as Big Data analysis or cosmological background wave investigations.

\section{Acknowledgements}
This material is based upon work supported by the Air Force Office of Scientific Research under
award number FA9550-22-1-0007.

\printbibliography

\end{document}